\documentclass[12pt,english]{article}
\usepackage[T1]{fontenc}
\usepackage[latin9]{inputenc}
\usepackage{color}
\usepackage{amsmath}
\usepackage{amssymb}
\usepackage{esint}
\usepackage{babel}

{\large{}{}{}{}{}{}{}{}}

\global\long\def\cbr#1{\left\{  #1\right\}  }
\global\long\def\rbr#1{\left(#1\right)}

\global\long\def\E{\mathbb{E}}
\global\long\def\P{\mathbb{P}}
\global\long\def\R{\mathbb{R}}

\global\long\def\dd#1{\textnormal{d}#1}

\global\long\def\ra{\rightarrow}
\global\long\def\TTV#1#2#3{\text{TV}^{#3}\!\rbr{#1,#2}}
\global\long\def\TTVemph#1#2#3{\emph{TV}^{#3}\!\rbr{#1,#2}}

\global\long\def\ns{\infty}

\newtheorem{theorem}{Theorem} 
 
\newtheorem{lem}[theorem]{Lemma}
\newtheorem{rem}[theorem]{Remark}

\newenvironment{proof}{\par\noindent{\bf Proof.}}{\par\rightline{$\blacksquare$}}

\title{Optimal uniform approximation of Lévy
processes on Banach spaces with finite variation
processes}

\author{Witold M. Bednorz\footnote{University of Warsaw}, Rafa{\l} M. \L ochowski\footnote{Warsaw School of Economics and the University of Warsaw}  and Rafa{\l} Martynek \footnote{University of Warsaw}}

\begin{document}
\maketitle

\begin{abstract}
For a general càdlàg Lévy process on a separable Banach
space $V$ we estimate values of $\inf_{Y\in{\cal A}_{X}}\E\left\{ \psi\left(\Vert X-Y\Vert_{\infty}\right)+\TTV Y{\left[0,T\right]}{}\right\} $,
where ${\cal A}_{X}$ is the family of processes on $V$ adapted to
the natural filtration of $X$, $\psi$ has polynomial growth and
$\TTV Y{\left[0,T\right]}{}$ denotes the total variation of the process
$Y$ on the interval $[0,T]$. Next, we apply obtained estimates in
three specific cases: Brownian motion with drift on $\R$, standard
Brownian motion on $\R^{d}$ and a symmetric $\alpha$-stable process
($\alpha\in(1,2)$) on $\R$.
\end{abstract}

\section{Introduction and formulation of the problem}

Let $X_{t},$ $t\ge0,$ be a càdlàg Lévy process  $X$ attaining its values
in a separable Banach space $V$ (i.e. a process with a.s. càdlàg
paths and independent and stationary increments), and let ${\cal A}_{X}$
be the family of $V$-valued processes $Y_{t},$ $t\ge0,$ adapted
to the natural filtration of $X.$ By $\left|\cdot\right|$ we denote
the norm in $V.$ For $T>0$ and two processes $Y,Z:\Omega\times[0,+\ns)\ra V$
we denote 
\[
\Vert Y-Z\Vert_{\infty,\left[0,T\right]}:=\sup_{0\leq t\leq T}\left|Y_{t}-Z_{t}\right|
\]
and 
\[
\TTV Y{\left[0,T\right]}{}:=\sup_{n}\sup_{0\le t_{0}<t_{1}<\cdots<t_{n}\le T}\sum_{i=1}^{n}\left|Y_{t_{i}}-Y_{t_{i-1}}\right|
\]
that is for $\omega \in \Omega$, $\TTV {Y(\omega)}{\left[0,T\right]}{}$ is the total variation of the trajectory $Y\rbr{\omega}$ on the interval $[0,T]$.

In this paper we deal with the following optimisation problem. Given
are $T>0$ and a non-decreasing function $\psi:\left[0,+\infty\right)\ra\left[0,+\infty\right)$
calculate (or estimate up to universal constants)

\begin{equation}
V_{X}\left(\psi\right):=\inf_{Y\in{\cal A}_{X}}\E\left\{ \psi\left(\Vert X-Y\Vert_{\infty,\left[0,T\right]}\right)+\TTV Y{\left[0,T\right]}{}\right\} .\label{eq:opt_pr1}
\end{equation}
To make the problem non-trivial we assume that $\E\left|X_{1}\right|<+\ns.$

{This type of optimization  problems {appear naturally in several situations. For example, in financial models with small proportional transaction costs where $X$ is the process representing optimal investment strategy on (frictionless) market without transaction costs, while  $Y$ is the approximation of $X$ and the total variation of $Y$ is proportional to the transactions costs of the implementation of the strategy $Y$}, see for example \cite{KallsenMuhle2013}. This type of optimization  problems have no unified,
algorithmic solution since the generator of the total variation functional
is not well defined. Moreover, we deal with very general Lévy processes
attaining their values in general Banach spaces. However, using well known results of the renewal theory, some
ad-hoc reasoning, results obtained for the functional called }{\emph{truncated
variation}}{ and assuming that $\psi$ grows no faster
than some polynomial, we will be able to estimate (\ref{eq:opt_pr1})
up to universal constants, depending on $\psi$ in terms
of the characteristics of the process $X.$ Together with the estimates
we will provide the construction of the process $Z$ uniformly approximating
$X,$ for which these estimates hold.}

From the triangle inequality we immediately get that if $\Vert X-Y\Vert_{\infty,\left[0,T\right]}\le c/2$
then for any $0\le s\le t\le T,$ $\left|Y_{t}-Y_{s}\right|\ge\max\left\{ \left|X_{t}-X_{s}\right|-c,0\right\} ,$
thus 
\begin{align}
 & \TTV Y{\left[0,T\right]}{} \nonumber \\
 & \ge\TTV X{\left[0,T\right]}c=\sup_{n}\sup_{0\le t_{0}<t_{1}<\cdots<t_{n}\le T}\sum_{i=1}^{n}\max\left\{ \left|X_{t_i}-X_{t_{i-1}}\right|-c,0\right\}. \label{eq:lower_bound_tv}
\end{align}
The quantity on the right side of (\ref{eq:lower_bound_tv}) is called
the truncated variation of $X.$ In the case when $V=\R,$ from the
results of \cite[Remark 15]{LochowskiGhomrasniMMAS:2015} it is possible
to prove that for any $c>0$ there exists a process $X^{c}\in{\cal A}_{X}$
such that $\Vert X-X^{c}\Vert_{\infty,\left[0,T\right]}\le c/2$ and
\[
\TTV X{\left[0,T\right]}c\le\TTV{X^{c}}{\left[0,T\right]}{}\le\TTV X{\left[0,T\right]}c+c,
\]
thus in the case $V=\R$ we have the estimate 
\begin{eqnarray*}
 &  & \inf_{c>0}\left\{ \psi\left(\frac{c}{2}\right)+\E\TTV X{\left[0,T\right]}c\right\} \\
 &  & \le\inf_{Y\in{\cal A}_{X}}\E\left\{ \psi\left(\Vert X-Y\Vert_{\infty,\left[0,T\right]}\right)+\TTV Y{\left[0,T\right]}{}\right\} \\
 &  & \le\inf_{c>0}\left\{ \psi\left(\frac{c}{2}\right)+\E\TTV X{\left[0,T\right]}c+ c\right\} ,
\end{eqnarray*}
which means that if {$\psi\left(x\right)$ grows no
faster than some polynomial (and no slower than some increasing linear function)
then} $\inf_{c>0}\left\{ \psi\left(c/2\right)+\E\TTV X{\left[0,T\right]}c\right\} $ and $V_{X}\left(\psi\right)$ 
are comparable up to universal
constants depending on $\psi$ only.

For a general Banach space-valued Lévy process, using similar construction
as in the proof of Theorem 1 from \cite{LochowskiJIA:2018},
we get that there exists a process $Y^{c}\in{\cal A}_{X}$ such that
$\Vert X-Y^{c}\Vert_{\infty,\left[0,T\right]}\le c/2$ and 
\begin{align}
\E\TTV X{\left[0,T\right]}c & \le\E\TTV{Y^{c}}{\left[0,T\right]}{}\label{optY}\\
 & \le\inf_{\lambda>1}\lambda\cdot\E\TTV X{\left[0,T\right]}{(\lambda-1)\cdot c/\left(2\lambda\right)}.\nonumber 
\end{align}
From this, assuming that there exists a constant $K_{\psi}$ such that for any $a\ge0,$ $\psi\left(2a\right)\le K_{\psi}\cdot\psi\left(a\right),$
we get 
\begin{align}
 & \inf_{c>0}\left\{ \psi\left(\frac{c}{2}\right)+\E\TTV X{\left[0,T\right]}c\right\} \nonumber \\
 & \le\inf_{Y\in{\cal A}_{X}}\E\left\{ \psi\left(\Vert X-Y\Vert_{\infty,\left[0,T\right]}\right)+\TTV Y{\left[0,T\right]}{}\right\} \nonumber \\
 & \le\inf_{c>0}\left\{ \psi\left(\frac{c}{2}\right)+\E\TTV{Y^{c}}{\left[0,T\right]}{}\right\} \label{eq:double_sided_opt_problem}\\
 & \le\inf_{c>0}\inf_{\lambda>1}\left\{ \psi\left(\frac{c}{2}\right)+\lambda\cdot\E\TTV X{\left[0,T\right]}{(\lambda-1)\cdot c/\left(2\lambda\right)}\right\} \nonumber \\
  & = \inf_{c>0}\inf_{\lambda>1}\left\{ \psi\left(\frac{4 c}{2}\right)+\lambda\cdot\E\TTV X{\left[0,T\right]}{(\lambda-1)\cdot 4 c/\left(2\lambda\right)}\right\} \nonumber \\
 & \le\inf_{c>0}\left\{ \psi\left(\frac{4c}{2}\right)+2\cdot\E\TTV X{\left[0,T\right]}{4c/4}\right\} \nonumber \\
 & \le\max\left(K_{\psi}^{2},2\right)\inf_{c>0}\left\{ \psi\left(\frac{c}{2}\right)+\E\TTV X{\left[0,T\right]}c\right\} \nonumber 
\end{align}
thus again we see that both quantities: $\inf_{c>0}\left\{ \psi\left({c}/{2}\right)+\E\TTV X{\left[0,T\right]}c\right\} $
and $V_{X}\left(\psi\right)$ are comparable up to universal
constants (depending on $\psi$ only). Since $X$ has càdlàg trajectories,
the construction of the process $Y^{c}$ appearing in (\ref{optY})
and (\ref{eq:double_sided_opt_problem}) simplifies to the following
one. First, we define stopping times $\tau_{0}^{c}=0$ and for $n=1,2,\ldots$
\begin{equation}
\tau_{n}^{c}=\begin{cases}
\inf\left\{ t>\tau_{n-1}^{c}:\left|X_{\tau_{n-1}^{c}}-X_{t}\right|\ge\frac{c}{2}\right\}  & \mbox{if }\tau_{n-1}^{c}<+\ns;\\
+\ns & \mbox{if }\tau_{n-1}^{c}=+\ns
\end{cases}\label{stop_times}
\end{equation}
and then we define 
\begin{equation}
Y_{t}^{c}=\sum_{n=0}^{+\ns}X_{\tau_{n}^{c}}\mathbf{1}_{\left[\tau_{n}^{c};\tau_{n+1}^{c}\right)}\left(t\right).\label{eq:yc_def}
\end{equation}
To avoid technical problems with these definitions we apply the
convention that $\inf\emptyset=+\ns$, $X_{\infty}=X_{0}$ and that
$\left[+\ns;+\ns\right)=\emptyset.$ 

\begin{rem} The construction in the proof of \cite[Theorem 1]{LochowskiJIA:2018}
rather uses times $\tilde{\tau}_{n}^{c}$ defined in the following
way 
\[
\tilde{\tau}_{n}^{c}=\begin{cases}
\inf\left\{ t>\tilde{\tau}_{n-1}^{c}:\left|X_{\tilde{\tau}_{n-1}^{c}}-X_{t}\right|>\frac{c}{2}\right\}  & \mbox{if }\tilde{\tau}_{n-1}^{c}<+\ns;\\
+\ns & \mbox{if }\tilde{\tau}_{n-1}^{c}=+\ns,
\end{cases}
\]
which may be not stopping times, but it is straightforward to verify
that for the times defined by (\ref{stop_times}) and $Y^{c}$ defined
by (\ref{eq:yc_def}) the estimates (\ref{eq:double_sided_opt_problem})
hold as well (see the proof of \cite[Theorem 1]{LochowskiJIA:2018}). 
\end{rem}
In what follows, we will use the presented construction to obtain
more straightforward estimates of $V_{X}\left(\psi\right)$
in terms of the characteristics of the process $X.$

{This paper is organised as follows. In the next section
we prove useful estimates of }$\mathbb{E}\TTV{Y^{c}}{\left[0,T\right]}{}$,
where $Y^{c}$ is the process defined by equation (\ref{eq:yc_def}),
and then prove two universal estimates of $V_{X}\left(\psi\right)$
(Theorem \ref{Theorem-1.-} and Theorem \ref{Theorem-2.-}) expressed
in terms of simpler functionals of $X$. In the last, third section,
we apply obtained estimates in three specific cases, namely when: (1) $X$ is a Brownian motion
with drift on $\R$, (2) $X$ is a standard Brownian motion on $\R^{d}$ and (3) $X$ is a symmetric
$\alpha$-stable process ($\alpha\in(1,2)$) on $\R$.

\section{Estimation of $\E\protect\TTV{Y^{c}}{\left[0,T\right]}{}$
and $V_{X}\left(\psi\right).$}

First, using the strong Markov property and the independence of the
increments of the process $X,$ we will estimate $\mathbb{E}\TTV{Y^{c}}{\left[0,T\right]}{}$,
where $Y^{c}$ is the process defined by equation (\ref{eq:yc_def}).

For $t>0$ let us define few auxiliary quantities
\[
\sigma^{c}(t):=\min\cbr{k:\tau_{k}^{c}>t} = \sum_{n=1}^{+\ns} \mathbf{1}_{\left\{ \tau_{n-1}^{c}\le t\right\} },
\]
\[
U^{c}(t):=\E\sigma^{c}(t) = \sum_{n=1}^{+\ns} \P\rbr{\tau_{n-1}^{c}\le t}  = \sum_{n=1}^{+\ns} \P\rbr{\sigma^{c}(t) \ge n}
\]
and 
\[
f^{c}(t):=\E \rbr{\left|X_{\tau_{1}^{c}}-X_{0}\right|\mathbf{1}_{\left\{ \tau_{1}^{c}\le t\right\} }}
\]
where stopping times $\tau_{n}^{c},$ $n=0,1,\ldots,$
are defined by formula (\ref{stop_times}).
\begin{lem}\label{lema0} For the process $Y^{c}$ defined by equation
(\ref{eq:yc_def}) the following inequalities hold:
\[
\frac{1}{2}U^{c}(T)f^{c}(T)\le\E\TTVemph{Y^{c}}{[0,T]}{}\le U^{c}(T)f^{c}(T).
\]
\end{lem}

\begin{proof}
Let us first notice that 
\begin{align*}
\TTV{Y^{c}}{\left[0,T\right]}{} & =\sum_{n=1}^{+\ns}\left|X_{\tau_{n}^{c}}-X_{\tau_{n-1}^{c}}\right|\mathbf{1}_{\left\{ \tau_{n}^{c}\le T\right\} },
\end{align*}
where stopping times $\tau_{n}^{c},$ $n=0,1,\ldots,$ are defined
by formula (\ref{stop_times}). 
From this, using independence of increments of $A$ and the strong Markov property we get an upper bound for
$\E\TTV{Y^{c}}{\left[0,T\right]}{},$ which reads
\begin{align*}
 & \E\TTV{Y^{c}}{\left[0,T\right]}{}=\E\sum_{n=1}^{+\ns}\left|X_{\tau_{n}^{c}}-X_{\tau_{n-1}^{c}}\right|\mathbf{1}_{\left\{ \tau_{n}^{c}\le T\right\} }\nonumber \\
 & \le\E\sum_{n=1}^{+\ns}\left|X_{\tau_{n}^{c}}-X_{\tau_{n-1}^{c}}\right|\mathbf{1}_{\left\{ \tau_{n}^{c}-\tau_{n-1}^{c}\le T\right\} }\mathbf{1}_{\left\{ \tau_{n-1}^{c}\le T\right\} } \nonumber \\
 & =\sum_{n=1}^{+\ns}\E\left(\left|X_{\tau_{n}^{c}}-X_{\tau_{n-1}^{c}}\right|\mathbf{1}_{\left\{ \tau_{n}^{c}-\tau_{n-1}^{c}\le T\right\} }\right)\E\mathbf{1}_{\left\{ \tau_{n-1}^{c}\le T\right\} } \\
 & =\sum_{n=1}^{+\ns}\E\left(\left|X_{\tau_{1}^{c}}-X_{0}\right|\mathbf{1}_{\left\{ \tau_{1}^{c}\le T\right\} }\right)\E\mathbf{1}_{\left\{ \tau_{n-1}^{c}\le T\right\} } \nonumber \\
 & = U^{c}(T)f^{c}(T).\nonumber 
\end{align*}
To bound $\E\TTV{Y^{c}}{\left[0,T\right]}{}$ from below we write 
\begin{align}
\TTV{Y^{c}}{\left[0,T\right]}{} & =\sum_{n=1}^{\sigma^{c}(T)-1}\left|X_{\tau_{n}^{c}}-X_{\tau_{n-1}^{c}}\right|\nonumber \\
 & = \sum_{n=1}^{\sigma^{c}(T)-1}\left|X_{\tau_{n}^{c}}-X_{\tau_{n-1}^{c}}\right|\mathbf{1}_{\left\{ \tau_{n}^{c}-\tau_{n-1}^{c}\le T\right\} }\mathbf{1}_{\left\{ \tau_{n-1}^{c}\le T\right\} }.\label{eq:first}
\end{align}
We will use the notion of \emph{stochastic domination}. We say that a real random variable
$Q$ stochastically dominates a real random variable $P$ if for any $x\in\R,$ $\P\left(Q\ge x\right)\ge\P\left(R\ge x\right).$
We denote this by $Q\succeq P$. We have that (even if $\P\left(\tau_{\sigma^{c}(T)}^{c}=+\ns\right)>0$, applying
the convention that $X_{\infty}=X_{0}$)
\begin{align}
\TTV{Y^{c}}{\left[0,T\right]}{} & \succeq \left|X_{\tau_{\sigma^{c}(T)}^{c}}-X_{\tau_{\sigma^{c}(T)-1}^{c}}\right|\mathbf{1}_{\left\{ \tau_{\sigma^{c}(T)}^{c}-\tau_{\sigma^{c}(T)-1}^{c}\le T\right\} }\mathbf{1}_{\left\{ \tau_{\sigma^{c}(T)-1}^{c}\le T\right\} }.\label{eq:second}
\end{align}
Taking expectations of both sides of relations \eqref{eq:first} and \eqref{eq:second} and adding them we get
\begin{align*}
2 \cdot  \E \TTV{Y^{c}}{\left[0,T\right]}{} & \ge \E \sum_{n=1}^{\sigma^{c}(T)}\left|X_{\tau_{n}^{c}}-X_{\tau_{n-1}^{c}}\right|\mathbf{1}_{\left\{ \tau_{n}^{c}-\tau_{n-1}^{c}\le T\right\} }\mathbf{1}_{\left\{ \tau_{n-1}^{c}\le T\right\} }\nonumber \\
 & = \E \sum_{n=1}^{+\ns }\left|X_{\tau_{n}^{c}}-X_{\tau_{n-1}^{c}}\right|\mathbf{1}_{\left\{ \tau_{n}^{c}-\tau_{n-1}^{c}\le T\right\} }\mathbf{1}_{\left\{ \tau_{n-1}^{c}\le T\right\} } \nonumber \\
 & = U^{c}(T)f^{c}(T),
 \end{align*}
where we used the fact that for $n > \sigma^{c}(T)$, $\mathbf{1}_{\left\{ \tau_{n-1}^{c}\le T\right\} } \equiv 0$  
\end{proof}
Let us denote $ \tau^c = \tau_1^c$. 
\begin{rem} \label{ren_fun}
The function $U^{c}(T)$ is example of a \emph{renewal function}, a well known object in the \emph{renewal theory}. \emph{Elementary renewal theorem} states that 
\[
\lim_{T \ra +\ns} \frac{U^{c}(T)}{T} = \frac{1}{\E \tau^c},
\]
where in the case $\E \tau^c = +\ns$ we set $1/\E\tau^c = 0$.
\end{rem}
\begin{rem} \label{eric}
We have the following estimates which are special case of results obtained by Erickson in \cite{Erickson1973}:
\[
\frac{t}{m^c(t)} \le U^{c}(T) \le \frac{2t}{m^c(t)},
\]
where $m^c(t) =\E \rbr{ \tau^c \wedge t } = \int_0^t \P \rbr{\tau^c > s} \dd s$.
This gives in the case $\E \tau^c = +\ns$ the proper order of growth of the renewal function.
\end{rem}
Sometimes (and  this is often the case when one deals with L\'evy processes) it is easier to deal with the Laplace transform of $\tau^c$ than with the function $U^{c}(T)$.  
\begin{lem} \label{lema}
For the process $Y^{c}$ defined by equation (\ref{eq:yc_def})
the following inequalities hold
\begin{equation}
\E\TTVemph{Y^{c}}{\left[0,T\right]}{}\le2\frac{f^c(T)}{1-\E2^{-\tau^{c}/T}}\label{eq:rel1}
\end{equation}
and 
\begin{equation}
\E\TTVemph{Y^{c}}{\left[0,T\right]}{}\ge\frac{1}{4}\frac{f^c(T)}{1-\E2^{-\tau^{c}/T}}.\label{eq:rel2}
\end{equation}
\end{lem}
\begin{proof}
Both estimates follow from Lemma \ref{lema0} and elementary estimates of $U^{c}(T)$. 
The estimate from above follows from the estimate
\begin{equation} \label{eq:estim_above}
U^{c}(T) = \sum_{n=1}^{+\ns} \E \mathbf{1}_{\cbr{\tau_{n-1}^{c}\le T}} \le \sum_{n=1}^{+\ns} \E 2^{1 - \tau_{n-1}^{c}/T}
\end{equation}
which is the consequence of the elementary estimate $\mathbf{1}_{\left\{ x\le T\right\} }\le 2^{1-x/T}$
valid for any $x\in\R.$ Further, we have 
\begin{align}
\sum_{n=1}^{+\ns}\E 2^{1 - \tau_{n-1}^{c}/T} & =2\sum_{n=1}^{+\ns}\E2^{- \tau_{n-1}^{c}/T}\nonumber \\
 & =2\sum_{n=1}^{+\ns}\left(\E 2^{-\tau^{c}/T}\right)^{n-1}\nonumber \\
 & =\frac{2}{1-\E 2^{-\tau^{c}/T}}.\label{eq:power_series}
\end{align}
From (\ref{eq:estim_above}) and (\ref{eq:power_series}) we get estimate
(\ref{eq:rel1}).

To bound $\E\TTV{Y^{c}}{\left[0,T\right]}{}$ from below for $t>0$ we define $\sigma^{c,0}(t)=0$ and for $k=1,2,\ldots,$ such that $\tau_{\sigma^{c,k-1}(t)}^c<+\ns$ let $\sigma^{c,k}(t)$ be the smallest integer such that $\tau_{\sigma^{c,k}(t)}^{c}-\tau_{\sigma^{c,k-1}(t)}^{c}>t$ (we naturally have $\sigma^{c,1}(t) = \sigma^{c}(t)$ and also have $\tau_{\sigma^{c,k}(t)-1}^{c}-\tau_{{\sigma^{c,k-1}(t)}}^{c}\le t$). For $k=1,2,\ldots,$ such that  $\tau_{\sigma^{c,k-1}(t)}^c = +\ns$ we set $\sigma^{c,k}(t) = \sigma^{c,k-1}(t)+1$. This yields that $\tau_{\sigma^{c,k}(t)}^{c} \ge k\cdot t$  
and for $k=0,1,2,\ldots$ we have 
\begin{align}
2^{-k} U^{c}(T) & \ge \E  \rbr{2^{-\tau_{\sigma^{c,k}(T)}^{c}/T}  \sum_{n= \sigma^{c,k}(T)+1}^{\sigma^{c,k+1}(T)} 1}.\label{eq:first1}
\end{align}
Summing estimates \eqref{eq:first1} over $k=0,1,2,\ldots$ we have 
\begin{align}
2 U^{c}(T)   &  = \sum_{k=0}^{+\ns} 2^{-k} U^{c}(T) \ge  \sum_{k=0}^{+\ns}  \E  \rbr{2^{-\tau_{\sigma^{c,k}(T)}^{c}/T}  \sum_{n= \sigma^{c,k}(T)+1}^{\sigma^{c,k+1}(T)} 1} \nonumber \\
& \ge \sum_{k=0}^{+\ns}  \E  \rbr{ 2^{-\tau_{\sigma^{c,k}(T)}^{c}/T} \sum_{n= \sigma^{c,k}(T)+1}^{\sigma^{c,k+1}(T)}  2^{-\rbr{\tau_{n-1}^{c}/T - \tau_{\sigma^{c,k}(T)}^{c}/T} }} \nonumber \\
& = \sum_{k=0}^{+\ns}  \E  \sum_{n= \sigma^{c,k}(T)+1}^{\sigma^{c,k+1}(T)}  2^{-\tau_{n-1}^{c}/T }  = \sum_{n=1}^{+\ns}  \E  2^{-\tau_{n-1}^{c}/T }\nonumber \\
& =\sum_{n=1}^{+\ns}\left(\E 2^{-\tau^{c}/T}\right)^{n-1} 
 =\frac{1}{1-\E 2^{-\tau^{c}/T}}.\label{eq:secondd}
\end{align}
Lemma \ref{lema0} and \eqref{eq:secondd} yield the estimate from below \eqref{eq:rel2}.
\end{proof}

Now, using Lemma \ref{lema} and estimates \eqref{eq:double_sided_opt_problem} we obtain the following result.
\begin{theorem}\label{Theorem-1.-} Let $X_{t},$ $t\ge0,$ be a Lévy process
on a separable Banach space $V$ with the norm $\left|\cdot\right|$
and let ${\cal A}_{X}$ be the class of processes adapted to the natural
filtration of $X.$ Let $\psi:\left[0,+\ns\right)\ra\left[0,+\ns\right)$
be a non-decreasing function such that for $a\ge0,$ $\psi\left(2a\right)\le K_{\psi}\cdot\psi\left(a\right).$
For any $T>0$ the following estimates hold: 
\begin{align}
V_{X}\left(\psi\right) & :=\inf_{Y\in{\cal A}_{X}}\E\left\{ \psi\left(\Vert X-Y\Vert_{\infty,\left[0,T\right]}\right)+\TTVemph Y{\left[0,T\right]}{}\right\} \nonumber \\
 & \le\inf_{c>0}\left\{ \psi\left(\frac{c}{2}\right)+ 2\frac{\E\left(\left|X_{\tau^{c}}-X_{0}\right|\mathbf{1}_{\left\{ \tau^{c}\le T\right\} }\right)}{1-\E 2^{-\tau^{c}/T}}\right\} \label{eq:thm1_est_above}\\
 & \le 2\inf_{c>0}\left\{ \psi\left(\frac{c}{2}\right)+\frac{\E\left(\left|X_{\tau^{c}}-X_{0}\right|\mathbf{1}_{\left\{ \tau^{c}\le T\right\} }\right)}{1-\E 2^{-\tau^{c}/T}}\right\} \nonumber
\end{align}
and 
\begin{align}
V_{X}\left(\psi\right) & =\inf_{Y\in{\cal A}_{X}}\E\left\{ \psi\left(\Vert X-Y\Vert_{\infty,\left[0,T\right]}\right)+\TTVemph Y{\left[0,T\right]}{}\right\} \nonumber \\
 & \ge\frac{1}{\max\left(K_{\psi}^{2},2\right)}\inf_{c>0}\left\{ \psi\left(\frac{c}{2}\right)+ \frac{1}{4}\frac{\E\left(\left|X_{\tau^{c}}-X_{0}\right|\mathbf{1}_{\left\{ \tau^{c}\le T\right\} }\right)}{1-\E 2^{-\tau^{c}/T}}\right\} \label{eq:thm1_est_below}\\
 & \ge\frac{1}{4\max\left(K_{\psi}^{2},2\right)}\inf_{c>0}\left\{ \psi\left(\frac{c}{2}\right)+\frac{\E\left(\left|X_{\tau^{c}}-X_{0}\right|\mathbf{1}_{\left\{ \tau^{c}\le T\right\} }\right)}{1-\E 2^{-\tau^{c}/T}}\right\}, \nonumber
\end{align}
where $\tau^{c}=\inf\left\{ t>0:\left|X_{t}-X_{0}\right|\ge c/2\right\} .$
\end{theorem}
In what follows we will estimate $\E\left|X_{\tau^{c}}-X_{0}\right|\mathbf{1}_{\left\{ \tau^{c}\le T\right\} }$
to obtain the following theorem.

\begin{theorem} \label{Theorem-2.-} Let $X_{t},$ $t\ge0,$ be a Lévy process
on a separable Banach space $V$ with the norm $\left|\cdot\right|$
and let ${\cal A}_{X}$ be the class of processes adapted to the natural
filtration of $X.$ Let $\psi:\left[0,+\ns\right)\ra\left[0,+\ns\right)$
be a non-decreasing function such that for $a\ge0,$ $\psi\left(2a\right)\le K_{\psi}\cdot\psi\left(a\right).$
For any $T,>0$ the following estimates hold: 
\begin{align}
 & V_{X}\left(\psi\right):=\inf_{Y\in{\cal A}_{X}}\E\left\{ \psi\left(\Vert X-Y\Vert_{\infty,\left[0,T\right]}\right)+\TTV Y{\left[0,T\right]}{}\right\} \label{eq:thm2_estim_above}\\
 & \le\inf_{c>0}\left\{ \psi\left(\frac{c}{2}\right)+3\frac{c\cdot\P\left(\tau^{c}\le T\right)}{1-\E2^{-\tau^{c}/T}}+ \frac{4}{\ln 2} T\int_{\left(c,+\ns\right]}y\cdot\Pi\left(\dd y\right)\right\} \nonumber 
\end{align}
and 
\begin{align}
 & V_{X}\left(\psi\right)=\inf_{Y\in{\cal A}_{X}}\E\left\{ \psi\left(\Vert X-Y\Vert_{\infty,\left[0,T\right]}\right)+\TTV Y{\left[0,T\right]}{}\right\} \label{eq:thm2_estim_below}\\
 & \ge\frac{1}{\max\left(K_{\psi}^{2},2\right)}\inf_{c>0}\left\{ \psi\left(\frac{c}{2}\right)+\frac{1}{16}\frac{c\cdot\P\left(\tau^{c}\le T\right)}{1-\E2^{-\tau^{c}/T}}+\frac{1}{32\ln 2}T\int_{\left(c,+\ns\right]}y\cdot\Pi\left(\dd y\right)\right\} \nonumber 
\end{align}
where $\tau^{c}=\inf\left\{ t>0:\left|X_{t}-X_{0}\right|\ge c/2\right\} $
and $\Pi$ is the image of the Lévy measure of the process $X$ under
the transformation $x\mapsto\left|x\right|$.
\end{theorem} 
\begin{proof} Let $\Delta X_{\tau^{c}}=X_{\tau^{c}}-X_{\tau^{c}-}$ denotes
the jump of the process at the moment $\tau^{c}$ and let us notice
that by the triangle inequality and the definition of $\tau^{c},$
for $\tau^{c}<+\ns$ we have 
\[
\left|X_{\tau^{c}}-X_{0}\right|\ge\left|\Delta X_{\tau^{c}}\right|-\left|X_{\tau^{c}-}-X_{0}\right|\ge\left|\Delta X_{\tau^{c}}\right|-c/2.
\]
Thus for $\tau^{c}<+\ns$ it follows that 
\begin{equation*} 
\left|\Delta X_{\tau^{c}}\right|\le c/2+\left|X_{\tau^{c}}-X_{0}\right|\le2\left|X_{\tau^{c}}-X_{0}\right|
\end{equation*}
and we have 
\begin{equation} 
\left|X_{\tau^{c}}-X_{0}\right|\ge\frac{1}{2}\left|\Delta X_{\tau^{c}}\right|.\label{eq:jump_est}
\end{equation}
Let now $\mu$ be the joint law of $\left(\left|\Delta X_{\tau^{c}}\right|,\tau^{c}\right)$.
For $y\in\left(c,+\ns\right)$ and $t\in\left(0,+\ns\right)$ one
has 
\[
\dd{\mu}\left(y,t\right)=\P\left(\sup_{0\le s<t}\left|X_{s}-X_{0}\right| < c/2\right)\dd{\Pi\left(y\right)}\dd t,
\]
where $\dd t$ denotes the Lebesgue measure. This observation follows
from the fact that for $y\in\left(c,+\ns\right)$ and $t\in\left(0,+\ns\right)$
the event 
\[
\left\{ \left|\Delta X_{\tau^{c}}\right|\in\left[y,y+\dd y\right),\tau^{c}\in\left[t,t+\dd t\right)\right\} 
\]
is equal the intersection of two independent events 
\[
\left\{ \sup_{0\le s<t}\left|X_{s}-X_{0}\right| < c/2\right\} \text{ and }\left\{ \left|X_{t+\dd t} - X_{t-} \right|\in\left[y,y+\dd y\right)\right\} 
\]
which follows from (\ref{eq:jump_est}) and the Lévy-Ito decomposition
(see \cite{Applebaum:2007}). Now, using \eqref{eq:jump_est} we easily estimate 
\begin{align}
\E\left|X_{\tau^{c}}-X_{0}\right|\mathbf{1}_{\left\{ \tau^{c}\le T\right\} } & \ge\E\left|X_{\tau^{c}}-X_{0}\right|\mathbf{1}_{\left\{ \left|\Delta X_{\tau^{c}}\right|>c\right\} }\mathbf{1}_{\left\{ \tau^{c}\le T\right\} }\nonumber \\
 & \ge\frac{1}{2}\E\left|\Delta X_{\tau^{c}}\right|\mathbf{1}_{\left\{ \left|\Delta X_{\tau^{c}}\right|>c\right\} }\mathbf{1}_{\left\{ \tau^{c}\le T\right\} }\nonumber \\
 & =\frac{1}{2}\int_{\left(c,+\ns\right]\times\left(0,T\right]}y\cdot\dd{\mu}\left(y,t\right)\nonumber \\
 & =\frac{1}{2}\int_{\left(c,+\ns\right]\times\left(0,T\right]}y\cdot\P\left(\sup_{0\le s<t}\left|X_{s}-X_{0}\right| < c/2\right)\dd{\Pi\left(y\right)}\dd t\nonumber \\
 & =\frac{1}{2}\int_{\left(c,+\ns\right]}y\cdot\Pi\left(\dd y\right)\int_{\left(0,T\right]}\P\left(\sup_{0\le s<t}\left|X_{s}-X_{0}\right| < c/2\right)\dd t\nonumber \\
 & =\frac{1}{2}\int_{\left(c,+\ns\right]}y\cdot\Pi\left(\dd y\right)\int_{\left(0,T\right]}\P\left(\tau^{c}\ge t\right)\dd t. \label{eq:one}
\end{align}
We naturally also have 
\begin{equation}
\E\left|X_{\tau^{c}}-X_{0}\right|\mathbf{1}_{\left\{ \tau^{c}\le T\right\} }\ge\frac{1}{2}c\cdot\P\left(\tau^{c}\le T\right).\label{eq:two}
\end{equation}
From (\ref{eq:one}) and (\ref{eq:two}) we get 
\begin{equation}
\E\left|X_{\tau^{c}}-X_{0}\right|\mathbf{1}_{\left\{ \tau^{c}\le T\right\} }\ge\frac{1}{4}c\cdot\P\left(\tau^{c}\le T\right)+\frac{1}{4}\int_{\left(c,+\ns\right]}y\cdot\Pi\left(\dd y\right)\int_{\left(0,T\right]}\P\left(\tau^{c}\ge t\right)\dd t.\label{eq:four}
\end{equation}
On the other hand, by the definition of $\tau^{c},$ for $\tau^{c}<+\ns$
we have 
\begin{align*}
\left|X_{\tau^{c}}-X_{0}\right| & \le\left|X_{\tau^{c}-}-X_{0}\right|+\left|\Delta X_{\tau^{c}}\right|\\
 & \le\frac{1}{2}c+\left|\Delta X_{\tau^{c}}\right|
\end{align*}
from which we get the estimate 
\begin{align}
&  \E\left|X_{\tau^{c}}-X_{0}\right|\mathbf{1}_{\left\{ \tau^{c}\le T\right\} } \le\frac{1}{2}c \cdot \E\mathbf{1}_{\left\{ \tau^{c}\le T\right\} }+\E\left|\Delta X_{\tau^{c}}\right|\mathbf{1}_{\left\{ \tau^{c}\le T\right\} }\nonumber \\
 & =\frac{1}{2}c\cdot\P\left(\tau^{c}\le T\right)+\E\left|\Delta X_{\tau^{c}}\right|\mathbf{1}_{\left\{ \left|\Delta X_{\tau^{c}}\right|\le c\right\} }\mathbf{1}_{\left\{ \tau^{c}\le T\right\} }\nonumber  \\
 & \quad+\E\left|\Delta X_{\tau^{c}}\right|\mathbf{1}_{\left\{ \left|\Delta X_{\tau^{c}}\right|>c\right\} }\mathbf{1}_{\left\{ \tau^{c}\le T\right\} }\nonumber \\
 & \le\frac{3}{2}c\cdot\P\left(\tau^{c}\le T\right)+\int_{\left(c,+\ns\right]}y\cdot\Pi\left(\dd y\right)\int_{\left(0,T\right]}\P\left(\tau^{c}\ge t\right)\dd t.\label{eq:three}
\end{align}
To deal with the integral $\int_{\left(0,T\right]}\P\left(\tau^{c}\ge t\right)\dd t$
let us notice that the following estimates hold: 
\[
\int_{\left(0,T\right]}\P\left(\tau^{c}\ge t\right)\dd t\le 2\int_{0}^{+\ns}2^{-t/T}\P\left(\tau^{c}\ge t\right)\dd t
\]
and 
\begin{align*}
\int_{0}^{+\ns}2^{-t/T}\P\left(\tau^{c}\ge t\right)\dd t & =\sum_{k=1}^{+\ns}\int_{(k-1)T}^{kT}2^{-t/T}\P\left(\tau^{c}\ge t\right)\dd t\\
 & \le\sum_{k=1}^{+\ns}\int_{(k-1)T}^{kT}2^{-(k-1)T/T}\P\left(\tau^{c}\ge t-(k-1)T\right)\dd t\\
 & =\sum_{k=1}^{+\ns}2^{-(k-1)}\int_{0}^{T}\P\left(\tau^{c}\ge t\right)\dd t\\
 & =2\int_{\left(0,T\right]}\P\left(\tau^{c}\ge t\right)\dd t.
\end{align*}
Thus, we have the double-sided estimate 
\begin{align}
\frac{1}{2}\int_{0}^{+\ns}2^{-t/T}\P\left(\tau^{c}\ge t\right)\dd t & \le\int_{\left(0,T\right]}\P\left(\tau^{c}\ge t\right)\dd t\nonumber \\
 & \le 2\int_{0}^{+\ns}2^{-t/T}\P\left(\tau^{c}\ge t\right)\dd t.\label{eq:double_sid_thm2}
\end{align}
Finally, let us notice that (by integration by parts) 
\begin{align}
\int_{0}^{+\ns}2^{-t/T}\P\left(\tau^{c}\ge t\right)\dd t & =\frac{T}{\ln 2}-\frac{T}{\ln 2}\int_{0}^{+\ns}2^{-t/T}\P\left(\tau^{c}\in\dd t\right)\nonumber \\
 & =\frac{T}{\ln 2}\left(1-\E2^{-\tau^{c}/T}\right).\label{eq:int_parts_thm2}
\end{align}

Now, from (\ref{eq:thm1_est_above}), (\ref{eq:three}), (\ref{eq:double_sid_thm2})
and (\ref{eq:int_parts_thm2}) we get (\ref{eq:thm2_estim_above})
while from (\ref{eq:thm1_est_below}), (\ref{eq:four}), (\ref{eq:double_sid_thm2})
and (\ref{eq:int_parts_thm2}) we get (\ref{eq:thm2_estim_below}).
\end{proof}

\section{Examples}

In this section we will apply the obtained estimates in three special
cases. In the first case the process $X$ will be a real-valued Brownian
motion with drift, in the second case it will be a standard Brownian
motion on $\R^{d}$, $d=2,3,\ldots,$ and in the third case it will be a real valued,
symmetric $\alpha$-stable process with $\alpha \in (1,2)$.

\subsection{Estimates of $V_{X}\left(\psi\right)$ in the case when $X$
is a Brownian motion with drift}

Let now $B$ be a standard Brownian motion starting from $0$ and
$X_{t}=B_{t}+\mu t$ be a (real-valued) Brownian motion with drift
$\mu.$ From Theorem \ref{Theorem-1.-} it follows that in order to estimate $V_{X}\left(\psi \right)$
it is sufficient to estimate (up to universal constants) the quantity
$\E\left(\left|X_{\tau^{c}}\right|\mathbf{1}_{\left\{ \tau^{c}\le T\right\} }\right)/\left(1-\E\exp\left(-\tau^{c}/T\right)\right).$
From the continuity of Brownian paths we immediately get that $\left|X_{\tau^{c}}\right|=c/2$
and 
\[
\E\left(\left|X_{\tau^{c}}\right|\mathbf{1}_{\left\{ \tau^{c}\le T\right\} }\right)=\frac{c}{2}\mathbb{P}\left(\tau^{c}\le T\right)=\frac{c}{2}\mathbb{P}\left(\sup_{0\le t\le T}\left|B_{t}+\mu t\right|\ge\frac{c}{2}\right).
\]
Now let us consider two cases.

1. $\frac{c}{2}\le\sqrt{T}+\left|\mu\right|T.$ Let $\tilde{B}=\text{sign\ensuremath{\left(\mu\right)}}B,$
where $\text{sign}{\left(\mu\right)} = -1 $ if $\mu <0$, $\text{sign}{\left(\mu\right)} = 1 $ if $\mu \ge 0$. 
We get 
\begin{align*}
\mathbb{P}\left(\sup_{0\le t\le T}\left|B_{t}+\mu t\right|\ge\frac{c}{2}\right) & \ge\mathbb{P}\left(\left|B_{T}+\mu T\right|\ge\sqrt{T}+\left|\mu\right|T\right)\\
 & =\mathbb{P}\left(\left|\text{sign\ensuremath{\left(\mu\right)}}B_{T}+\text{sign\ensuremath{\left(\mu\right)}}\mu T\right|\ge\sqrt{T}+\left|\mu\right|T\right)\\
 & \ge\mathbb{P}\left(\text{sign\ensuremath{\left(\mu\right)}}B_{T}+\left|\mu\right|T\ge\sqrt{T}+\left|\mu\right|T\right)\\
 & =\mathbb{P}\left(\tilde{B}_{T}\ge\sqrt{T}\right)\ge1-\Phi\left(1\right)> \frac{1}{7},
\end{align*}
where $\Phi\left(x\right)=\left(2\pi\right)^{-1/2}\int_{-\infty}^{x}e^{-t^{2}/2}\mathrm{d}t$
is the cumulative probability function of a standard normal variable.

2. $\frac{c}{2}>\sqrt{T}+\left|\mu\right|T.$ In this case we get
the following lower bound 
\begin{align*}
\mathbb{P}\left(\sup_{0\le t\le T}\left|B_{t}+\mu t\right|\ge\frac{c}{2}\right) & \ge\mathbb{P}\left(\left|B_{T}+\mu T\right|\ge\frac{c}{2}\right)\\
 & =\mathbb{P}\left(\left|\text{sign\ensuremath{\left(\mu\right)}}B_{T}+\text{sign\ensuremath{\left(\mu\right)}}\mu T\right|\ge\frac{c}{2}\right)\\
 & \ge\mathbb{P}\left(\text{sign\ensuremath{\left(\mu\right)}}B_{T}+\left|\mu\right|T\ge\frac{c}{2}\right)\\
 & =\mathbb{P}\left(\tilde{B}_{T}\ge\frac{c}{2}-\left|\mu\right|T\right)=1-\Phi\left(\frac{\frac{c}{2}-\left|\mu\right|T}{\sqrt{T}}\right).
\end{align*}
On the other hand, in both cases we have 
\begin{align*}
\mathbb{P}\left(\sup_{0\le t\le T}\left|B_{t}+\mu t\right|\ge\frac{c}{2}\right) & \le\mathbb{P}\left(\sup_{0\le t\le T}\left|B_{s}\right|+\left|\mu\right|T\ge\frac{c}{2}\right)\\
 & \le2\mathbb{P}\left(\sup_{0\le t\le T}B_{t}\ge\frac{c}{2}-\left|\mu\right|T\right)\\
 & =4\mathbb{P}\left(B_{T}\ge\frac{c}{2}-\left|\mu\right|T\right)\\
 & =4\left(1-\Phi\left(\frac{\frac{c}{2}-\left|\mu\right|T}{\sqrt{T}}\right)\right).
\end{align*}
Thus, in both cases 
\[
\mathbb{P}\left(\sup_{0\le t\le T}\left|B_{t}+\mu t\right|\ge\frac{c}{2}\right)=\kappa_{1}\left(1-\Phi\left(\frac{\frac{c}{2}-\left|\mu\right|T}{\sqrt{T}}\right)\right),
\]
where $\kappa_{1}\in\left[\frac{1}{7},4\right].$

Now we turn to look at $1-\E 2^{-\tau^{c}/T}=\E\left(1-e^{-\ln 2 \cdot \tau^{c} /T}\right).$
Let $\tau$ be an exponentially distributed random variable, independent
from $B,$ with the cumulative probability function $\mathbb{P}\left(\tau<t\right)=\left(1-e^{-\ln 2 \cdot t/T}\right)\mathbf{1}_{\{t>0\}}.$
By formula 1.15.2 on p. 270 in \cite{BorodinSalminen:2002} we get
\begin{align*}
\E\left(1-2^{-\tau^{c}/T}\right) & =\mathbb{P}\left(\tau<\tau^{c}\right)\\
 & =\mathbb{P}\left(\inf_{0\le s\le\tau}\left(B_{s}+\mu s\right)>-\frac{c}{2},\sup_{0\le s\le\tau}\left(B_{s}+\mu s\right)<\frac{c}{2}\right)\\
 & =1-\frac{\left(e^{-\mu c/2}+e^{\mu c/2}\right)\sinh\left(\frac{c}{2}\sqrt{\frac{2 \ln 2}{T}+\mu^{2}}\right)}{\sinh\left(c\sqrt{\frac{2 \ln 2}{T}+\mu^{2}}\right)} \\
 & = 1 - \frac{\cosh\rbr{\frac{\mu c}{2}} }{\cosh\left(\frac{c}{2}\sqrt{\frac{2 \ln 2}{T}+\mu^{2}}\right)}.
\end{align*}

Finally, from Theorem \ref{Theorem-1.-} we get that 
\[
V_{X}\left(\psi\right)=\kappa_{2}\inf_{c>0}\left\{ \psi\left(\frac{c}{2}\right)+\frac{c}{2}\frac{1-\Phi\left(\left(\frac{c}{2}-\left|\mu\right|T\right)/\sqrt{T}\right)}{1 - \frac{\cosh\rbr{\frac{\mu c}{2}} }{\cosh\left(\frac{c}{2}\sqrt{\frac{2 \ln 2}{T}+\mu^{2}}\right)}}\right\} ,
\]
where $\kappa_{2}\in\left[\frac{1}{28\max\left(K_{\psi}^{2},2\right)},8\right].$

\subsection{Estimates of $V_{X}\left(\psi\right)$ in the case when $X$ is a
$d$-dimensional standard Brownian motion ($d\ge2$)}

Let $B^{\left(1\right)},B^{\left(2\right)},\ldots,B^{\left(d\right)}$
be $d$ independent ($d\ge2$), standard, one-dimensional Brownian
motions, starting from $0$, and let $X=\left(B^{\left(1\right)},B^{\left(2\right)},\ldots,B^{\left(d\right)}\right)$.

Again, by Theorem \ref{Theorem-1.-} it is sufficient to estimate
the ratio $$\frac{\E\left(\left|X_{\tau^{c}}\right|\mathbf{1}_{\left\{ \tau^{c}\le T\right\} }\right)}{1-\E2^{-\tau^{c}/T}}$$
and again, by the continuity of $X,$ we get 
\[
\E\left(\left|X_{\tau^{c}}-X_{0}\right|\mathbf{1}_{\left\{ \tau^{c}\le T\right\} }\right)=\frac{c}{2}\mathbb{P}\left(\tau^{c}\le T\right).
\]
Moreover, recall that the process $R$ defined by 
\[
R=\sqrt{\left(B^{\left(1\right)}\right)^{2}+\left(B^{\left(2\right)}\right)^{2}+\ldots+\left(B^{\left(d\right)}\right)^{2}}
\]
is called $d$-dimensional Bessel process or a Bessel process of order
$d$ or a Bessel process with index $\nu=d/2-1$. 

Using results of G. Serafin \cite{Serafin:2017}, we will obtain estimates
of $V_{X}\left(\psi\right)$ which are universal up to a constant
depending on $\psi$ and $\nu$ (but not $T$). Using \cite[Corollary 3.4]{Serafin:2017}
and scaling properties of the standard Brownian motion for $t>0$
we get 
\[
\P\rbr{\tau^{c}\in\dd t}=\kappa\left(t,\nu\right)\frac{4}{c^{2}}\rbr{1+\frac{c^{2}}{4t}}^{\nu+2}\exp\rbr{-\frac{1}{8}\frac{c^{2}}{t}-2j_{\nu,1}^{2}\frac{t}{c^{2}}}\dd t,
\]
where $\kappa(t,\nu)\in[\kappa_{3}\left(\nu\right),\kappa_{4}\left(\nu\right)]$
and $\kappa_{4}\left(\nu\right)>\kappa_{3}(\nu)>0$ are constants
depending on $\nu$ only, and $j_{\nu,1}$ denotes the smallest positive
zero of the the Bessel function $J_{\nu}$ of the first kind which
is defined as
\[
J_{\nu}\left(y\right)=\left(\frac{y}{2}\right)^{\nu}\sum_{m=0}^{+\infty}\frac{\left(-1\right)^{m}}{m!\Gamma\left(m+\nu+1\right)}\left(\frac{y}{2}\right)^{2m}.
\]
This gives the following estimates.

1. $c^{2}\le T.$ In this case 
\begin{align*}
\P\rbr{\tau^{c}\le T} & =\int_{0}^{T}\P\rbr{\tau^{c}\in\dd t}\ge\int_{c^{2}/2}^{c^{2}}\P\rbr{\tau^{c}\in\dd t}\\
 & \ge\kappa_{3}(\nu)\frac{4}{c^{2}}\int_{c^{2}/2}^{c^{2}}\rbr{1+\frac{c^{2}}{4t}}^{\nu+2}\exp\rbr{-\frac{1}{8}\frac{c^{2}}{t}-2j_{\nu,1}^{2}\frac{t}{c^{2}}}\dd t\\
 & \ge\kappa_{3}(\nu)\frac{4}{c^{2}}\int_{c^{2}/2}^{c^{2}}\exp\rbr{-\frac{1}{4}-2j_{\nu,1}^{2}}\dd t\\
 & =\kappa_{5}(\nu).
\end{align*}

2. $c^{2}>T.$ In this case 
\begin{align*}
\P\rbr{\tau^{c}\le T} & =\int_{0}^{T}\P\rbr{\tau^{c}\in\dd t}\\
 & \ge\kappa_{3}(\nu)\frac{4}{c^{2}}\int_{0}^{T}\rbr{1+\frac{c^{2}}{4t}}^{\nu+2}\exp\rbr{-\frac{1}{8}\frac{c^{2}}{t}-2j_{\nu,1}^{2}\frac{t}{c^{2}}}\dd t\\
 & \ge\kappa_{3}(\nu)\frac{4}{c^{2}}\int_{0}^{T}\rbr{\frac{c^{2}}{4t}}^{\nu+2}\exp\rbr{-\frac{1}{8}\frac{c^{2}}{t}-2j_{\nu,1}^{2}}\dd t\\
 & \ge\kappa_{6}(\nu)\frac{4}{c^{2}}\int_{0}^{T}\rbr{\frac{c^{2}}{4t}}^{\nu+2}\exp\rbr{-\frac{1}{8}\frac{c^{2}}{t}}\dd t\\
 & =\kappa_{7}(\nu)\Gamma\rbr{\nu+1,\frac{c^{2}}{8T}},
\end{align*}
where for $a,y>0$, $\Gamma\rbr{a,y}$ denotes the incomplete gamma
function, 
\[
\Gamma\rbr{a,y}=\int_{y}^{+\infty}x^{a-1}e^{-x}\dd x.
\]
Similarly we can obtain a bound from above:
\begin{align*}
\P\rbr{\tau^{c}\le T} & =\int_{0}^{T}\P\rbr{\tau^{c}\in\dd t}\\
 & \le\kappa_{4}(\nu)\frac{4}{c^{2}}\int_{0}^{T}\rbr{1+\frac{c^{2}}{4t}}^{\nu+2}\exp\rbr{-\frac{1}{8}\frac{c^{2}}{t}-2j_{\nu,1}^{2}\frac{t}{c^{2}}}\dd t\\
 & \le\kappa_{4}(\nu)\frac{4}{c^{2}}\int_{0}^{T}\rbr{4\frac{c^{2}}{4t}+\frac{c^{2}}{4t}}^{\nu+2}\exp\rbr{-\frac{1}{8}\frac{c^{2}}{t}}\dd t\\
 & \le\kappa_{8}(\nu)\frac{4}{c^{2}}\int_{0}^{T}\rbr{\frac{c^{2}}{4t}}^{\nu+2}\exp\rbr{-\frac{1}{8}\frac{c^{2}}{t}}\dd t\\
 & =\kappa_{9}(\nu)\Gamma\rbr{\nu+1,\frac{c^{2}}{8T}}.
\end{align*}

We notice that in both ($c^2 \le T$ and $c^2 >T$) cases
\begin{equation}
\P\rbr{\tau^{c}\le T}=\kappa_{10}\Gamma\rbr{\nu+1,\frac{1}{8}\max\rbr{1,\frac{c^{2}}{T}}},\label{eq:Bessel2-1}
\end{equation}
where $\kappa_{10}\in[\kappa_{11}\left(\nu\right),\kappa_{12}\left(\nu\right)]$
and $\kappa_{12}\left(\nu\right)>\kappa_{11}(\nu)>0$ are constants
depending on $\nu$ only. 

Next, the term $1-\E2^{-\tau^{c}/T} = 1-\E\exp\left(-(\ln 2 /T) \cdot \tau^{c}\right)$ is given by an explicit
formula. 
For a Bessel process with index $\nu$, starting from $x$ and $\lambda \ge 0$:  
\[
1-\E^{x}\exp\left(-\lambda \cdot \tau^{c}\right)=1-\left(\frac{c}{2}\right)^{\nu}\frac{x^{-\nu}I_{\nu}\left(x\sqrt{2\lambda }\right)}{I_{\nu}\left(c\sqrt{\lambda/2}\right)}
\]
(see \cite[formula 1.1.2 on p. 373]{BorodinSalminen:2002}). 
Here $I_{\nu}$ denotes the modified Bessel function 
\[
I_{\nu}\left(y\right)=\left(\frac{y}{2}\right)^{\nu}\sum_{m=0}^{+\infty}\frac{1}{m!\Gamma\left(m+\nu+1\right)}\left(\frac{y}{2}\right)^{2m}
\]
and in our case (for $x=0$) we get $x^{-\nu}I_{\nu}\left(x\sqrt{2\lambda}\right)=\left(\lambda/2\right)^{\nu/2}/\Gamma\left(\nu+1\right),$
hence, substituting $\lambda = \ln 2 /T$, 
\begin{align}
1-\E2^{-\tau^{c}/T} & =1-\frac{\left(c^2\ln2/(8T)\right)^{\nu/2}}{\Gamma\left(\nu+1\right)I_{\nu}\left(c\sqrt{\ln2/(2T)}\right)}.\label{eq:Bessel1}
\end{align}

(\ref{eq:Bessel2-1}) and (\ref{eq:Bessel1}) together with Theorem
\ref{Theorem-1.-} allow to estimate $V_{X}\left(\psi\right)$ up
to a constants depending on $\nu$ and $\psi$ only: 
\begin{align*}
 & V_{X}\left(\psi\right)\le2\inf_{c>0}\left\{ \psi\left(\frac{c}{2}\right)+\frac{c}{2}\frac{\kappa_{12}\left(\nu\right)\Gamma\rbr{\nu+1,\frac{1}{8}\max\rbr{1,\frac{c^{2}}{T}}}}{1-\frac{\left(c^2\ln2/(8T)\right)^{\nu/2}}{\Gamma\left(\nu+1\right)I_{\nu}\left(c\sqrt{\ln2/(2T)}\right)}}\right\} 
\end{align*}
and 
\[
V_{X}\left(\psi\right)\ge\frac{1}{4\max\left(K_{\psi}^{2},2\right)}\inf_{c>0}\left\{ \psi\left(\frac{c}{2}\right)+\frac{c}{2}\frac{\kappa_{11}\left(\nu\right)\Gamma\rbr{\nu+1,\frac{1}{8}\max\rbr{1,\frac{c^{2}}{T}}}}{1-\frac{\left(c^2\ln2/(8T)\right)^{\nu/2}}{\Gamma\left(\nu+1\right)I_{\nu}\left(c\sqrt{\ln2/(2T)}\right)}}\right\} .
\]

\begin{rem}
By formula 1.1.4 on p. 373 in \cite{BorodinSalminen:2002},
which is due to J. T. Kent \cite{Kent:1980}, we get exact formula
for $\mathbb{P}\left(\tau^{c}\le T\right)$:
\begin{align}
\mathbb{P}\left(\tau^{c}\le T\right) & =1-\frac{2^{1-\nu}}{\Gamma\left(\nu+1\right)}\sum_{k=1}^{+\infty}\frac{j_{\nu,k}^{\nu-1}}{J_{\nu+1}\left(j_{\nu,k}\right)}e^{-2j_{\nu,k}^{2}T/c^{2}},\label{eq:Bessel}
\end{align}
where $j_{\nu,1} < j_{\nu,2} < \ldots$ denote consecutive positive zeros of $J_{\nu}$.
Unfortunately, the series in formula (\ref{eq:Bessel}) is convenient
when dealing with larger times, and for $T\ge2\mathbb{E}\tau^{c}=2c^{2}/d$
we naturally have $\mathbb{P}\left(\tau^{c}\le T\right)\ge1/2.$ Unfortunately,
for smaller $T$s this sum is oscillating and in that case the cancellations
between the terms really matter in the context of asymptotic behaviour. 
\end{rem}

\subsection{Estimates of $V_{X}\left(\psi\right)$ in the case when $X$
is a symmetric, real-valued, strictly $\alpha$-stable motion ($\alpha\in(1,2)$)}

Let now $X_{t}$, $t\ge0$, be a symmetric, real-valued, strictly
$\alpha$-stable motion ($\alpha\in(1,2)$) such that $X_{0}=0$.
$X$ has the following scaling property: for $t\ge0$ and $a>0$ 
\begin{equation}
X_{a t}\stackrel{law}=a^{1/\alpha}X_{t}.\label{eq:scaling}
\end{equation}
To fix our attention to the process of a given magnitude, we will
assume that $X_{1}$ has the following characteristic function 
\[
\E\exp\rbr{i \xi X_{1}}=\exp\rbr{\int_{-\ns}^{\ns}e^{i \xi x}-1-i \xi x\frac{\dd x}{\left|x\right|^{\alpha+1}}}=e^{-\sigma_{\alpha}\left|\xi \right|^{\alpha}},
\]
where $\xi \in \R$, $\sigma_{\alpha}=2\Gamma\left(-\alpha\right)\cos\frac{\left(2-\alpha\right)\pi}{2}$.
Let $\beta_{\alpha}$ be such that 
\[
\P\rbr{\left|X_{1}\right|\ge\beta_{\alpha}}=\frac{1}{3 e^5}.
\]

To estimate $V_{X}\left(\psi\right)$ we will apply Theorem \ref{Theorem-2.-}. 

First we need to estimate $1-\E 2^{-\tau^{c}/T}=\frac{\ln 2}{T}\int_{0}^{+\ns}2^{-t/T}\P\left(\tau^{c}\ge t\right)\dd t$.
Let us define the function $\left(0,+\ns\right)\ni c\mapsto u\left(c\right)\in\left(0,+\ns\right)$
such that 
\[
\P\rbr{\left|X_{u\rbr c}\right|\ge\frac{c}{2}}=\frac{1}{3 e^5}.
\]
 By scaling property of $X$ it is equivalent to 
\[
\P\rbr{\left|X_{1}\right|\ge \left(u\rbr c\right)^{-1/\alpha}\frac{c}{2}}=\frac{1}{3 e^5}=\P\rbr{\left|X_{1}\right|\ge\beta_{\alpha}}.
\]
Thus 
\begin{equation}
u\left(c\right)=\frac{c^{\alpha}}{\rbr{2\beta_{\alpha}}^{\alpha}}.\label{eq:uc}
\end{equation}
Next, by symmetry and strong Markov property of $X$, we have the
estimate 
\begin{align}
\frac{\P\rbr{\left|X_{u\rbr c}\right|\ge\frac{c}{2}}}{\P\rbr{\tau^{c}\le u\rbr c}} & =\P\rbr{\left|X_{u\rbr c}\right|\ge\frac{c}{2}|\tau^{c}\le u\rbr c}\nonumber \\
 & \ge\P\rbr{\text{sign}\left(X_{\tau^{c}}\right)\left(X_{u\rbr c}-X_{\tau^{c}}\right)\ge0|\tau^{c}\le u\rbr c}\nonumber \\
 & \ge\frac{1}{2}\label{eq:Levy}
\end{align}
(recall that $\left|X_{\tau^c}\right|\ge\frac{c}{2}$) from which it follows $\P\rbr{\tau^{c}\le u\rbr c}\le2\P\rbr{\left|X_{u\rbr c}\right|\ge\frac{c}{2}}$
and 
\begin{align}
\P\rbr{\tau^{c}>u\left(c\right)} & =1-\P\left(\tau^{c}\le u\rbr c\right)\nonumber \\
 & \ge1-2\P\rbr{\left|X_{u\rbr c}\right|\ge\frac{c}{2}}\nonumber \\
 & \ge1-2\frac{1}{3e^5}>\frac{1}{2}.\label{eq:estim_below-1}
\end{align}
On the other hand, using the independence of the increments and scaling
properties of $X$, for $k=1,2,\ldots$, we estimate 
\begin{align}
\P\rbr{\tau^{c}>2^{\alpha}k\cdot u\left(c\right)} & =\P\rbr{\sup_{0\le s\le2^{\alpha}k\cdot u\left(c\right)}\left|X_{s}\right|<\frac{c}{2}}\nonumber \\
 & \le\P\rbr{\sup_{0\le s<t\le2^{\alpha}k\cdot u\left(c\right)}\left|X_{s}-X_{t}\right|<c}\nonumber \\
 & \le\P\rbr{\sup_{2^{\alpha}\left(j-1\right)\cdot u\left(c\right)\le s<t\le2^{\alpha}j\cdot u\left(c\right)}\left|X_{s}-X_{t}\right|<c\text{ for }j=1,2,\ldots,k}\nonumber \\
 & =\prod_{j=1}^{k}\P\rbr{\sup_{2^{\alpha}\left(j-1\right)\cdot u\left(c\right)\le s<t\le2^{\alpha}j\cdot u\left(c\right)}\left|X_{s}-X_{t}\right|<c}\nonumber \\
 & =\left(\P\rbr{\sup_{0\le s<t\le2^{\alpha}u\left(c\right)}\left|X_{s}-X_{t}\right|<c}\right)^{k}\nonumber \\
 & \le\left(\P\rbr{\sup_{0\le s\le2^{\alpha}u\left(c\right)}\left|X_{s}\right|<c}\right)^{k}=\left(\P\rbr{\sup_{0\le s\le u\left(c\right)}2\left|X_{s}\right|<c}\right)^{k}\nonumber \\
 & =\left(1-\P\rbr{\sup_{0\le s\le u\left(c\right)}\left|X_{s}\right|\ge\frac{c}{2}}\right)^{k}\le\left(1-\P\rbr{\left|X_{u\left(c\right)}\right|\ge\frac{c}{2}}\right)^{k}\nonumber \\
 & =\rbr{1-\frac{1}{3e^5}}^{k}.\label{eq:estim_above-1}
\end{align}
From (\ref{eq:estim_below-1}) we estimate 
\begin{align}
1-\E2^{-\tau^{c}/T} & =\frac{\ln 2}{T}\int_{0}^{+\ns}2^{-t/T}\P\left(\tau^{c}\ge t\right)\dd t\nonumber \\
 & \ge\frac{\ln 2}{T}\int_{0}^{u(c)}2^{-t/T}\frac{1}{2}\dd t\nonumber \\
 & =\frac{1}{2} \left(1-2^{-u(c)/T} \right) \label{eq:1skladnik}
\end{align}
and from (\ref{eq:estim_above-1}) we estimate 
\begin{align}
1-\E2^{-\tau^{c}/T} & =\frac{\ln 2}{T}\int_{0}^{+\ns}2^{-t/T}\P\left(\tau^{c}\ge t\right)\dd t\nonumber \\
 & =\frac{\ln 2}{T} \sum_{k=1}^{+\ns}\int_{2^{\alpha}\left(k-1\right)\cdot u\left(c\right)}^{2^{\alpha}k\cdot u\left(c\right)}2^{-t/T}\P\left(\tau^{c}\ge t\right)\dd t\nonumber \\
 & \le \frac{\ln 2}{T} \sum_{k=1}^{+\ns}\int_{2^{\alpha}\left(k-1\right)\cdot u\left(c\right)}^{2^{\alpha}k\cdot u\left(c\right)}2^{-t/T}\rbr{1-\frac{1}{3e^5}}^{k-1}\dd t\nonumber \\
 & =\frac{1}{1-\rbr{1-\frac{1}{3e^5}}2^{-2^{\alpha}u\left(c\right)}}\left(1-e^{-2^{\alpha}u\left(c\right)/T}\right)\nonumber \\
 & \le 3e^5 \cdot2^{\alpha}\left(1-2^{-u\left(c\right)/T}\right)\le12 e^5\left(1-2^{-u\left(c\right)/T}\right).\label{eq:2skladnik}
\end{align}
The last but one inequality follows from the estimates: $1-\rbr{1-\frac{1}{3e^5}}2^{-2^{\alpha}u\left(c\right)} \ge\frac{1}{3e^5}$
and $1-2^{-2^{\alpha}u\left(c\right)/T}\ge2^{\alpha}\left(1-2^{-u\left(c\right)/T}\right)$
which is the consequence of the concavity of the function $x\mapsto1-2^{-x}$:
\[
\frac{1}{2^{\alpha}}\left(1-2^{-2^{\alpha}u\left(c\right)/T}\right)+\left(1-\frac{1}{2^{\alpha}}\right)\left(1-2^{-0}\right)\le1-2^{-\frac{1}{2^{\alpha}}2^{\alpha}u\left(c\right)/T}.
\]

Next, we need to estimate $\P\rbr{\tau^{c}\le T}=\P\rbr{\sup_{0\le s\le T}\left|X_{s}\right|\ge\frac{c}{2}}.$
Using similar reasoning as in (\ref{eq:Levy}) we get:
\[
\P\rbr{\left|X_{T}\right|\ge\frac{c}{2}}\le\P\rbr{\tau^{c}\le T}\le2\P\rbr{\left|X_{T}\right|\ge\frac{c}{2}}.
\]
 This and scaling properties of $X$ yield:
\begin{equation}
\P\rbr{\left|X_{1}\right|\ge\frac{c}{2T^{1/\alpha}}}\le\P\rbr{\tau^{c}\le T}\le2\P\rbr{\left|X_{1}\right|\ge\frac{c}{2T^{1/\alpha}}}.\label{eq:3skladnik}
\end{equation}
Finally, using Theorem \ref{Theorem-2.-}, (\ref{eq:3skladnik}),
(\ref{eq:1skladnik}), (\ref{eq:uc}) and the fact that $\Pi\left(\dd y\right)=\left|y\right|^{-\alpha-1}\dd y$,
we obtain estimate from above:
\begin{equation} \label{alfa_above}
V_{X}\left(\psi\right)\le\inf_{c>0}\left\{ \psi\left(\frac{c}{2}\right)+12\frac{c\cdot\P\rbr{\left|X_{1}\right|\ge\frac{c}{2T^{1/\alpha}}}}{1-2^{-c^{\alpha}/\left(\rbr{2\beta_{\alpha}}^{\alpha}T\right)}}+\frac{4}{\ln 2} \frac{T}{\left(\alpha-1\right)c^{\alpha-1}}\right\} .
\end{equation}
Similarly, using Theorem \ref{Theorem-2.-}, (\ref{eq:3skladnik}),
(\ref{eq:2skladnik}), (\ref{eq:uc}) and the fact that $\Pi\left(\dd y\right)=\left|y\right|^{-\alpha-1}\dd y$,
we obtain estimate from below:
\begin{align}
 & V_{X}\left(\psi\right)\ge \label{alfa_below}\\
 & \frac{1}{\max\left(K_{\psi}^{2},2\right)}\inf_{c>0}\left\{ \psi\left(\frac{c}{2}\right)+\frac{1}{16^2}\frac{\P\rbr{\left|X_{1}\right|\ge\frac{c}{2T^{1/\alpha}}}}{1-2^{-c^{\alpha}/\left(\rbr{2\beta_{\alpha}}^{\alpha}T\right)}}+\frac{1}{32 \ln 2}\frac{T}{\left(\alpha-1\right)c^{\alpha-1}}\right\} . \nonumber
\end{align}

From the proof of \cite[Theorem 13]{BednorzLochMartynek:2018} it follows that 
$\beta_{\alpha} \ge \frac{2}{\sqrt{2-\alpha}}$. From the same \cite[Theorem 13]{BednorzLochMartynek:2018} (estimate (48) in \cite{BednorzLochMartynek:2018}) it follows that $\beta_{\alpha} \le \frac{\tilde{D}}{\sqrt{2-\alpha}}$ for some universal constant $\tilde{D}$. 
From this and concavity of the function $x \mapsto 1- 2^{-x}$ we obtain that 
\begin{equation} \label{12x}
d \rbr{1-2^{-c^{\alpha}/\left(\rbr{ \frac{1}{\sqrt{2-\alpha}}}^{\alpha}T\right)}} \le 1-2^{-c^{\alpha}/\left(\rbr{2\beta_{\alpha}}^{\alpha}T\right)} \le D \rbr{ 1-2^{-c^{\alpha}/\left(\rbr{ \frac{1}{\sqrt{2-\alpha}}}^{\alpha}T\right)}}
\end{equation}
for some universal constants $0<d\le D$.

Moreover, from \cite[Theorem 13]{BednorzLochMartynek:2018} it also follows that 
\begin{equation} \label{alfa1above}
\P\rbr{\left|X_{1}\right|\ge y} \le  C  \min\cbr{1, \max\cbr{ \frac{1}{(y/4)^{\alpha}}, e^{-(2-\alpha) (y/4)^2 } } }
\end{equation}
and
\begin{equation} \label{alfa1below}
\P\rbr{\left|X_{1}\right|\ge y} \ge  c \min\cbr{1,  \max\cbr{ \frac{1}{\rbr{\sqrt{8}y}^{\alpha}}, e^{-(2-\alpha) \rbr{\sqrt{8}y}^2 } }}
\end{equation}
for some universal constants $0<c\le C$.

Thus, from \eqref{alfa_above}-\eqref{alfa1below} we get that for some constant $L_{\psi}$, depending on $\psi$ only, we have
\begin{equation} 
V_{X}\left(\psi\right)\le L_{\psi} \inf_{c>0}\left\{ \psi\left(\frac{c}{2}\right)+\frac{c\cdot F \rbr{\frac{c}{T^{1/\alpha} }}}{1-2^{-c^{\alpha}{\rbr{2-\alpha}}^{\alpha/2}/T}}+\frac{T}{\left(\alpha-1\right)c^{\alpha-1}}\right\}
\end{equation}
and
\begin{equation} 
V_{X}\left(\psi\right)\ge \frac{1}{L_{\psi}} \inf_{c>0}\left\{ \psi\left(\frac{c}{2}\right)+\frac{c\cdot F \rbr{\frac{c}{T^{1/\alpha} }}}{1-2^{-c^{\alpha}{\rbr{2-\alpha}}^{\alpha/2}/T}}+\frac{T}{\left(\alpha-1\right)c^{\alpha-1}}\right\},
\end{equation}
 where 
 $$
 F(y) = \min\cbr{1,  \max\cbr{ {y}^{-\alpha}, e^{-(2-\alpha) y^2 } }}.
 $$

{\bf Acknowledgments.} Fruitful disscussion with René Schilling helped to
prove Theorem \ref{Theorem-2.-}, Jan Palczewski indicated alternative
statement of the optimisation problem, which is currently used in
this article. The research of all authors was funded by the National
Science Centre, Poland, under Grant No. 2016/21/B/ST1/0148. The research of W. M. B. and R. M. {\L }. was also funded by the National Science Centre, Poland, under Grant No. 2019/35/B/ST1/042 while the research  of R. M. was also funded by the National Science Centre, Poland, under Grant No. 2018/31/N/ST1/03982.

\bibliographystyle{plain}
\bibliography{/Users/rafallochowski/biblio/biblio}

\end{document}